\documentclass[12pt]{article}
\usepackage{amssymb,amsfonts,amsmath, psfrag,eepic,colordvi,epsfig}
\textheight by \topskip \numberwithin{equation}{section}
\newtheorem{theorem}{Theorem}[section]

\newtheorem{conjecture}[theorem]{Conjecture}

\newtheorem{lemma}[theorem]{Lemma}

\def\qed{\hfill \rule{4pt}{7pt}}

\begin{document}
\parskip 10pt

\begin{center}

{\Large\bf  Families of Sets with Intersecting Clusters }

\vskip 3mm

William Y.C. Chen$^1$ \\
Center for Combinatorics, LPMC-TJKLC \\
Nankai University, Tianjin 300071, P. R. China

\vskip 1mm

Jiuqiang Liu$^2$ \\
Department of Mathematics\\
Eastern Michigan University \\
Ypsilanti, MI 48197, USA

\vskip 1mm

Larry X.W. Wang$^3$ \\
Center for Combinatorics, LPMC-TJKLC \\
Nankai University, Tianjin 300071, P. R. China

\vskip 1mm

$^1$chen@nankai.edu.cn, $^2$jliu@emich.edu,
$^3$wxw@cfc.nankai.edu.cn
\end{center}

\begin{center}
 {\it In Memory of Professor Chao Ko}
\end{center}

\begin{abstract}
A family of $k$-subsets $A_1, A_2, \ldots, A_d$ on $[n]=\{1,
2,\ldots, n\}$ is called a $(d, c)$-cluster if the union $A_1\cup
A_2 \cup \cdots \cup A_d$ contains at most $ck$ elements with $c<d$.
Let $\mathcal{F}$ be a family of $k$-subsets of an $n$-element set.
We show that for $k \geq 2$ and $n \geq k+2$, if every $(k,
2)$-cluster of $\mathcal{F}$ is intersecting,
  then $\mathcal{F}$
contains no $(k-1)$-dimensional simplices. This leads to an
affirmative answer to Mubayi's conjecture for $d=k$  based on
Chv\'atal's simplex theorem.   We also show that for any $d$
satisfying $3 \leq d \leq k$ and $n \geq \frac{dk}{d-1}$, if every
$(d, {d+1\over 2})$-cluster is intersecting,  then
$|\mathcal{F}|\leq {{n-1} \choose {k-1}}$ with equality only when $
\mathcal{F}$ is a complete star. This result is an extension of both
Frankl's theorem and Mubayi's theorem.
\end{abstract}

\noindent {\bf Keywords:} Clusters of subsets, Chv\'atal's simplex
theorem, $d$-simplex, Erd\H{o}s-Ko-Rado Theorem
\\

\noindent {\bf AMS Classification}: 05D05.

\section{Introduction}

This paper is concerned with the study of families of subsets with
intersecting clusters. The first result is a proof of an important
case of a conjecture recently proposed by Mubayi \cite{mubayi-a} on
intersecting families with the aid of Chv\'atal's simplex theorem.
The second result is an extension of both Frankl's theorem and
Mubayi's theorem. It should be noted that we have used these two
theorems themselves as a starting point to prove this extension.

Let us review some notation and terminology. The set $\{1, 2, \dots,
n\}$ is usually denoted by $[n]$ and the  family of all $k$-subsets
of a finite set $X$ is denoted by  $X^{k}$ or ${{X} \choose {k}}$. A
family $\mathcal{F}$ of sets is said to be intersecting if every two
sets in $\mathcal{F}$ have a nonempty intersection. A family
$\mathcal{F}$ of sets in $X^{k}$ is called a complete star if $
\mathcal{F}$ consists of all $k$-subsets containing $x$ for some
$x\in X$.

The classical  Erd\H{o}s-Ko-Rado (EKR) theorem \cite{ekr} is stated
as follows.

\begin{theorem}[The EKR Theorem]  Let $n \geq 2k$ and let $\mathcal{F}
\subseteq {{[n]} \choose {k}}$ be an intersecting family, then
$|\mathcal{F}|\leq {{n-1} \choose {k-1}}$. Furthermore, for $n>2k$,
the equality holds only when $\mathcal{F}$ is a complete star.
\end{theorem}

The following generalization of the EKR theorem is due to Frankl
\cite{frankl}.

\begin{theorem}[Frankl] Let $k \geq 2$, $d \geq 2$, and $n
\geq dk/(d-1)$. Suppose that $\mathcal{F} \subseteq [n]^k$ such
that every $d$ sets of $\mathcal{F}$ have a nonempty intersection.
Then $|\mathcal{F}|\leq {{n-1} \choose {k-1}}$ with equality only
when $\mathcal{F}$ is a complete star.
\end{theorem}

The following  conjecture due to Erd\H{o}s  on triangle free
families implies Frankl's theorem for $d \geq 3$. Recall that a
$d$-dimensional simplex, or a $d$-simplex for short, is defined to
be a family of $d+1$ sets $A_1, A_2, \dots, A_{d+1}$ such that every
$d$ of them have a nonempty intersection, but $A_1 \cap A_2 \cap
\cdots \cap A_{d+1} = \emptyset$.  A $2$-dimensional simplex is
called a triangle. This conjecture has been proved by Mubayi and
Verstra\"ete \cite{mubayi-c}

\begin{conjecture}[Erd\H{o}s]
For $n \geq \frac{3k}{2}$, if $\mathcal{F} \subseteq [n]^k$
contains no triangle, then $|\mathcal{F}|\leq {{n-1} \choose
{k-1}}$ with equality only when $\mathcal{F}$ is a complete star.
\end{conjecture}

However, as generalization of Erd\H{o}s' conjecture,  Chv\'atal
\cite{chvatal} proposed the following conjecture which   remains
open in general case.

\begin{conjecture}[Chv\'atal's Simplex Conjecture]
Let $k \geq d+1 \geq 3$, $n \geq k(d+1)/d$, and $\mathcal{F}
\subseteq [n]^k$. If $\mathcal{F}$ contains no $d$-dimensional
simplex, then $|\mathcal{F}|\leq {{n-1} \choose {k-1}}$ with
equality only when $\mathcal{F}$ is a complete star.
\end{conjecture}

Chv\'atal \cite{chvatal} has shown that it is true for $d=k-1$,
which we call Chv\'atal's simplex theorem.

\begin{theorem}[Chv\'atal's Simplex Theorem]  \label{chvatal-t}
 For $n \geq k+2 \geq 5$, if
$\mathcal{F} \subseteq [n]^k$ contains no $(k-1)$-dimensional
simplices, then $|\mathcal{F}|\leq {{n-1} \choose {k-1}}$ with
equality only when $\mathcal{F}$ is a complete star.
\end{theorem}

Frankl and F\"uredi \cite{ff} have shown that Chv\'atal's conjecture
holds for sufficiently large $n$.

\begin{theorem}[Frankl and F\"uredi]
For $k \geq d+2 \geq 4$, there exists $n_0$ such that for $n >
n_0$, if $\mathcal{F} \subseteq [n]^k$ contains no $d$-dimensional
simplices, then $|\mathcal{F}|\leq {{n-1} \choose {k-1}}$ with
equality only when $\mathcal{F}$ is a complete star.
\end{theorem}

As will be seen, a recent conjecture proposed by Mubayi
\cite{mubayi-a} is  related to Chv\'atal's simplex theorem. Here we
introduce the terminology of  clusters of subsets. A family of
$k$-subsets $A_1, A_2, \ldots, A_d$ of $[n]$ is called a $(d,
c)$-cluster if $|A_1\cup A_2\cup \cdots \cup A_d| \leq ck$, where
$c<d$ is a constant that may depend on $d$. A cluster is said to be
intersecting if their intersection is nonempty.

\begin{conjecture}[Mubayi's Conjecture] Let $k \geq d \geq 3$ and $n \geq
dk/(d-1)$. Suppose that $\mathcal{F} \subseteq [n]^{k}$ such that
every $(d,2)$-cluster of $\mathcal{F}$ is intersecting£¬ i.e., for
any $A_1, A_2, \dots, A_d \in \mathcal{F}$, $|A_1 \cup A_2 \cup
\cdots \cup A_d| \leq 2k$ implies $A_1 \cap A_2 \cap \cdots \cap
A_d \neq \emptyset$. Then $|\mathcal{F}|\leq {{n-1} \choose
{k-1}}$ with equality only when $\mathcal{F}$ is a complete star.
\end{conjecture}

Mubayi \cite{mubayi-a} has shown that this conjecture holds for
$d=3$ (Theorem \ref{Mubayi-t}). He has also proved that his
conjecture holds for $d=4$ when $n$ is sufficiently large
\cite{mubayi-b}.

\begin{theorem}[Mubayi]\label{Mubayi-t}
Let $k \geq 3$ and $n \geq \frac{3k}{2}$. Suppose that
$\mathcal{F} \subseteq [n]^k$ is a family such that every $(3,
2)$-cluster $A_1, A_2, A_3 \in \mathcal{F}$ is intersecting, then
$|\mathcal{F}|\leq {{n-1} \choose {k-1}}$ with equality only when
$\mathcal{F}$ is a complete star.
\end{theorem}

In this paper, we  study the case $d=k$ of Mubayi's conjecture in
connection with Chv\'atal's simplex theorem. We show that in this
case the conditions for Mubayi's conjecture imply the nonexistence
of any $(k-1)$-dimensional simplex. Therefore, Chv\'atal's simplex
theorem leads to Mubayi's conjecture for $d=k$. As the main result
of this paper, we present a theorem on families of subsets with
intersecting clusters which can be viewed as an extension of both
Frankl's Theorem (Theorem 1.2) and Mubayi's Theorem  (Theorem 1.8).

\section{Families of Subsets with Intersecting Clusters}

In this section, we first consider a special case of Mubayi's
conjecture for $k=d$. We show that this case can be deduced from
Chv\'atal's simplex theorem (Theorem \ref{chvatal-t}). Then we study
families of $k$-subsets with intersecting $(d, {d+1\over
2})$-clusters and obtain a theorem as an extension of both Frankl's
theorem (Theorem 1.2) and Mubayi's theorem (Theorem 1.8). Our proof
is based on the EKR Theorem and Frankl's Theorem. We will also use a
similar strategy as in the proof of Mubayi's theorem
\cite{mubayi-a}.

\begin{theorem}
 Let $k \geq 3$ and $n \geq k+2$. Suppose that $\mathcal{F}
\subseteq [n]^{k}$ is a family of subsets of $[n]$ such that every
$(k, 2)$-cluster is intersecting.  Then $\mathcal{F}$ contains no
$(k-1)$-dimensional simplices.
\end{theorem}

\noindent {\it Proof.}   Suppose to the contrary that $A_1, A_2,
\dots, A_k \in \mathcal{F}$ form a $(k-1)$-dimensional simplex,
namely, every $k-1$ of them have a nonempty intersection but
\begin{equation}\label{a1k}
A_1 \cap A_2 \cap \cdots \cap A_k = \emptyset.
\end{equation}
It follows that two distinct families $\{A_{i_1}, A_{i_2}, \ldots,
A_{i_{k-1}}\}$ and $\{A_{j_1}, A_{j_2}, \ldots, A_{j_{k-1}}\}$
cannot have a common element, because the union of these two
families equals $\{A_1, A_2, \ldots, A_k\}$. Without loss of
generality, let
\[
i \in A_1\cap \cdots \cap A_{i-1}\cap A_{i+1}\cap \cdots \cap A_{k}.
\]
That is, $i$ belongs to every subset $A_j$ other than $A_i$. It
follows that that $\{1,\ldots ,i-1,i+1,\ldots k\}\subset A_i$. Since
$A_i$ is a $k$-subset, $A_i$ must contain an element in
$\{k+1,\ldots ,n\}$.  So we have
\[
|A_1\cup A_2\cup \cdots \cup A_k|\leq 2k.
\]
This means that $\{A_1,A_2\ldots ,A_k\}$ is a $(k, 2)$-cluster that
is not intersecting, contradicting to the assumption of the theorem.
So we conclude that $\mathcal{F}$ does not contain any
$(k-1)$-dimensional simplex. This completes the proof. \qed

The following theorem is the main result of this paper.

\begin{theorem}
Let $k \geq d \geq 3$ and $n \geq \frac{dk}{d-1}$. Suppose that
$\mathcal{F} \subseteq [n]^k$ is a family of subsets of $[n]$ such
that every $(d, {d+1\over 2})$-cluster is intersecting (i.e., for
any $A_1, A_2, \dots, A_d \in \mathcal{F}$,  $|A_1\cup A_2\cup
\cdots \cup A_d| \leq {d+1\over 2}k$ implies that
$\cap_{i=1}^{d}A_i \neq \emptyset$). Then $|\mathcal{F}|\leq
{{n-1} \choose {k-1}}$ with equality only when $ \mathcal{F}$ is a
complete star.
\end{theorem}

The next  lemma gives an upper bound on the number of edges in a
graph with intersecting clusters, and it  will be used in the
proof of Theorem 2.2.

\begin{lemma} \label{l-l-2}
Let $n > d \geq 3$. Suppose that $\mathcal{F} \subseteq [n]^2$ is
a family of $2$-subsets of $[n]$ such that every $(d, {d+1\over
2})$-cluster is intersecting. Then $|\mathcal{F}|\leq n-1 $ with
equality only when $\mathcal{F}$ is a complete star.
\end{lemma}

\noindent {\it Proof.} Since $\mathcal{F}$ is a family of
$2$-subsets, we may consider it as a graph $G$ with vertex set $
[n]$. The conditions in the lemma imply that any
 $d$ edges $A_1$,
$A_2$, $\dots$, $A_d$ of $G$ either intersect at a common vertex or
cover at least $d+2$ vertices (for $d=3$, $G$ does not contain any
triangle because every $(3,2)$-cluster is intersecting).

We  proceed by induction on $n$. For $n = d+1$, since any $d$ edges
cover at most $n=d+1$ vertices, any $d$ edges of $G$ must intersect
at a common vertex and thus form a star. This implies that
$|\mathcal{F}|=|E(G)| \leq d = n-1$ with equality only when
$\mathcal{F}$ (or $G$) is a complete star.

Assume that $n \geq d+2$ and  that the lemma holds for $n - 1$. We
first claim that $G$ must contain a vertex of degree one. Otherwise,
every vertex of $G$ has degree at least two which implies that for
every connected component $C$ of $G$ we have
\begin{equation}\label{ve}
|V(C)| \leq |E(C)|.
 \end{equation}
Let $C_1$, $C_2$, $\dots$, $C_m$ be the connected components of $G$
ordered by the condition
\[|E(C_1)| \geq |E(C_2)|\geq\cdots \geq |E(C_m)|.\]

We aim to find $d$ edges that form a non-intersecting $(d, {d+1\over
2})$-cluster to reach a contradiction. Let us consider two cases.

Case 1. $|C_1|\geq d$. Since $C_1$ is not a star, it contains a path
$P$ with three edges. Since $d\geq 3$, we can add $d-3$ edges to $P$
to obtained a connected subgraph $H$ of $C_1$. Let $A_1,A_2,\ldots
,A_d  $ be $d$ edges of $H$. Then we have
\[
|A_1\cup A_2\ldots \cup A_d|=|V(H)| \leq |E(H)| + 1=d+1.
\]
Since $H$ is not a star, we obtain $A_1\cap A_2\ldots \cap
A_d=\emptyset$.

Case 2. $|C_1|<d$. Let $r \geq 1$ be the integer such that
\[ b=\sum_{i=1}^r |E(C_i)| <d \quad \mbox{and} \quad
\sum_{i=1}^{r+1} |E(C_i)| \geq d.\] It is clear that $C_{r+1}$ has
at least $d-b$ edges. We now take any connected subgraph $H$ of
$C_{r+1}$ with $d-b$ edges. Since $H$ is connected, we have
\begin{equation} \label{ev} |E(H)| \geq |V(H)| - 1.
\end{equation}
Let $A_1, A_2 , \dots, A_{d}$ be the $d$ edges in $C_1, C_2,
\ldots, C_r, H$. From (\ref{ve}) and (\ref{ev}) it follows that
\begin{eqnarray*}
\lefteqn{ |A_1 \cup A_2 \cdots \cup A_d| }\\[3pt]
& = &  |V(C_1)|+ |V(C_2)|+ \cdots +|V(C_r)| + |V(H)| \\[3pt]
& \leq & |E(C_1)|+ |E(C_2)|+ \cdots +|E(C_r)| + |E(H)|+1\\[3pt]
 & = & d+1.
 \end{eqnarray*}
Noting that $C_1, C_2, \ldots, C_r$ and $H$ are disjoint,  we have
$A_1 \cap A_2 \cdots \cap A_d = \emptyset$.

In summary, we have reached the conclusion that $G$ has a vertex
with degree one. Let $v$ be a vertex of degree one in $G$ and let
$G'$ be the induced graph obtained from $G$ by deleting the vertex
$v$. Clearly, $G'$ is a graph with $n-1$ vertices in which every
$d$ edges $A_1$, $A_2$, $\dots$, $A_d$ either intersect at a
common vertex or cover at least $d+2$ vertices. By the inductive
hypothesis, we have $|E(G')| \leq n-2$ with equality only if $G'$
is a complete star. Hence
\[ |\mathcal{F}|=|E(G)|=|E(C)|+ 1 \leq n-1\]
with equality only if $\mathcal{F}$ (or $G$) is a complete star.
 \qed

The following lemma is an extension  of Lemma 3 of Mubayi
\cite{mubayi-a}. While the proof of Mubayi relies on the EKR
theorem, our proof is based on the above Lemma \ref{l-l-2} and
Frankl's theorem (Theorem 1.2). We will also use a similar framework
as in the proof of Mubayi's theorem \cite{mubayi-a}.

\begin{lemma}
Let $k \geq d \geq 2$, $t \geq 2$, and $2 \leq l \leq k$. Let
$S_1$, $S_2$, $\dots$, $S_t$ be pairwise disjoint $k$-subsets and
$X = S_1 \cup S_2 \cup \cdots \cup S_t$. Suppose that
$\mathcal{F}$ is a family of $l$-subsets of $X$ satisfying  the
conditions (1)  $S_i \in \mathcal{F}$ for all $i$ if $l = k$;
 (2) For every $A_1$, $A_2$, $\dots$, $A_{d}$ $\in \mathcal{F}$
and $ 1\leq i \leq t$, $A_1 \cap A_2 \cdots \cap A_{d} \cap S_i =
\emptyset$ implies $|A_1 \cup A_2 \cdots \cup A_{d}-S_i|>
\frac{dl}{2}$. Then we have $|\mathcal{F}|< {{tk-1} \choose
{l-1}}$.
\end{lemma}

\noindent {\it Proof.} For $d = 2$, the above lemma reduces to
Lemma 3 in \cite{mubayi-a}. So we  may assume that $d \geq 3$. Let
$n = |X| = tk$. We consider the following two cases.

\noindent Case 1. Assume $l = 2$. We claim that any $(d, {d+1\over
2})$-cluster of $\mathcal{F}$ is intersecting, namely, for any
$A_1$, $A_2$, $\dots$, $A_{d } \in \mathcal{F}$, we have either
$A_1 \cap A_2 \cap \cdots \cap A_{d } \neq \emptyset$ or $ | A_1
\cup A_2 \cup \cdots \cup A_{d } | \geq d+2$. To this end, we
assume that $A_1 \cap A_2 \cap \cdots \cap A_{d } = \emptyset$.
This gives $A_1 \cap A_2 \cap \cdots \cap A_{d }\cap S_{i} =
\emptyset$ for any $S_i$. Since $X=\cup S_i$ is the ground set of
$\mathcal{F}$, there exists $S_m$ such that $A_1 \cap S_{m} \neq
\emptyset$. As $A_1 \cap A_2 \cap \cdots \cap A_{d }\cap S_{m} =
\emptyset$ and $l=2$, in view of Condition 2 we get
\[|A_1 \cup A_2 \cup \cdots \cup A_{d }-S_{m}|> d . \]
Furthermore,  the condition $A_1 \cap S_{m} \neq \emptyset$ yields
\[ |A_1 \cup A_2 \cup \cdots \cup A_{d }|> d+1 .\]  So the
claim holds.

Since $d  \geq 3$, by Lemma 2.3,  we find that $|\mathcal{F}|\leq
n-1$, where $n=tk$. So it remains to show that it is impossible for
$|\mathcal{F}|$ to reach the upper bound $n-1$. Assume that
$|\mathcal{F}|= n-1$. Again, by Lemma 2.3, $\mathcal{F}$ must be a
complete star, namely, $ \mathcal{F}$ consists of all $2$-subsets of
$X$ for some $x$ in $X$. Without loss of generality, we may assume
that $x \in S_1$. Let $A_1$ be a $2$-subset from $\mathcal{F}$ such
that $A_1\subseteq S_1$. Since $d-1 \leq k$, we may choose $d-1$
$2$-subsets $A_2, A_3, \ldots, A_d$ such that $A_i \in \mathcal{F}$
and $A_i-x \subseteq S_2$ for $2 \leq i \leq d$. This implies that
\[
A_1 \cap A_2 \cap \cdots \cap A_{d } \cap S_2 = \emptyset\]
 and
\[ |(A_1 \cup A_2 \cup \cdots \cup A_{d }) - S_2| = 2 < d ,\]
contradicting Condition (2). Thus we have $|\mathcal{F}|< n-1 =
tk-1$. So the lemma is proved for $l=2$.

\noindent Case 2. Assume $l\geq 3$. So we have $k \geq l \geq 3$. We
use  induction on $t$.

We first consider the  case $t=2$, namely, $X=S_1\cup S_2$.  We will
show that $A_1 \cap A_2 \cap \cdots \cap A_{d } \neq \emptyset$ for
any $A_1$, $A_2$, $\dots$, $A_{d }$ $\in \mathcal{F}$. If this were
not true,  there would exist subsets $A_1, A_2, \ldots, A_d \in
\mathcal{F}$ for which
\begin{equation} \label{a-e}
A_1 \cap A_2 \cap \cdots \cap A_{d } = \emptyset.
\end{equation}
Let $A=A_1\cup A_2 \cup \cdots \cup A_d$. It is clear that $A$
contains at most $dl$ elements. Since $S_1$ and $S_2$ are
disjoint, so are $A\cap S_1$ and $A \cap S_2$. Therefore, either
$A\cap S_1$ or $A\cap S_2$ contains at most half of the elements
in $A$. We may assume without loss of generality that
\[ | A \cap S_1| \leq \frac{d l}{2} .\]
Note that (\ref{a-e}) implies $A_1 \cap A_2 \cap \cdots \cap A_{d }
\cap S_2 = \emptyset$. Since $X=S_1\cup S_2$, we get
\[ |A-S_2|=|A \cap S_1| \leq \frac{d l}{2},
\] contradicting  Condition (2).  Thus we deduce that
 $A_1 \cap A_2 \cap \cdots \cap
A_{d } \neq \emptyset$ for any $A_1$, $A_2$, $\dots$, $A_{d }$
$\in \mathcal{F}$. By Frankl's Theorem (Theorem 1.2) we obtain
\begin{equation} \label{frankl-u}
|\mathcal{F}|\leq {{2k-1} \choose {l-1}}.
\end{equation}
Next we prove that the equality in (\ref{frankl-u}) can never be
reached. Let us assume that \begin{equation} \label{frankl-u-2}
|\mathcal{F}|={{2k-1} \choose {l-1}}.
\end{equation}
Since $d\geq 3$, by Frankl's theorem, $\mathcal{F}$ is a complete
star, that is, $\mathcal{F}$ consists of all  $l$-subsets of $[2k]$
containing an element $x$ for some $x$ in $[2k]$. Without loss of
generality, we may assume that $x \in S_1$. Thus $\mathcal{F}$
contains every subset $A_i$ which is either of the form $B\cup
\{x\}$ for $B\in [S_1-x]^{l-1}$ or of the form $C\cup \{x\}$ for
$C\in [S_2]^{l-1}$. Since $d \leq k$ and $3 \leq l \leq k$, we have
\[d-1 \leq k \leq {{k} \choose {l-1}}.\]
Now we may choose $A_1 \in \mathcal{F}$ with $A_1 \subseteq S_1$ and
$d-1$ sets $A_2$, $A_3$, $\dots$, $A_{d }$ $\in \mathcal{F}$ with
$A_i -{x} \subseteq S_2$ for each $i \geq 2$.  Since $A_1\cap
S_2=\emptyset$,  $A_1 \cap A_2 \cdots \cap A_{d } \cap S_2 =
\emptyset$. Moreover, since $A_i -x \subseteq S_2$ for $i=2, 3,
\ldots, d$, we have
\[ |(A_1 \cup A_2 \cup \cdots \cup A_{d }) - S_2|
= |A_1| = l < \frac{d l}{2},\] contradicting Condition (2). It
follows that $|\mathcal{F}|< {{2k-1} \choose {l-1}}$ and hence the
lemma is valid for $t=2$.

Next suppose that $t \geq 3$ and the result holds for $t-1$. We
first show that there exists at most one set $S_m$ such that
\[ | \mathcal{F}\cap [S_{m}]^l| \geq  \frac{d}{2}.\]
Suppose, to the contrary, that there exist two sets, say $S_1$ and
$S_2$ , such that
\[|\mathcal{F}\cap [S_i]^l| \geq \frac{d}{2},\]
for $i=1,2$. Then we have
\[
|\mathcal{F}\cap [S_1]^l| + |\mathcal{F}\cap [S_2]^l| \geq d.
\]
 Since $|\mathcal{F}\cap [S_1]^l| \geq \frac{d}{2} \geq 1$ and $|\mathcal{F}\cap [S_2]^l| \geq \frac{d}{2} \geq 1$,
  we are able to choose $d$ sets  $A_1, A_2, \ldots, A_d$ from
 $(\mathcal{F}\cap [S_1]^l)  \cup  (\mathcal{F}\cap [S_2]^l)$
 such that  $A_1 \subseteq S_1$ and $A_2 \subseteq S_2$.
 Since $|(A_1
\cup A_2 \cup \cdots \cup A_{d })|\leq  d l$ and $S_1\cap
S_2=\emptyset$, we have either
\begin{equation} \label{a-1} |(A_1
\cup A_2 \cup \cdots \cup A_{d }) \cap S_1| \leq \frac{d
l}{2}\end{equation}
 or \begin{equation} \label{a-2}
  |(A_1 \cup A_2 \cup \cdots \cup A_{d }) \cap
S_2| \leq \frac{d l}{2}.\end{equation} Without loss of generality,
assuming that (\ref{a-1}) is valid. We see that
\[ |(A_1 \cup A_2 \cup \cdots \cup A_{d }) - S_2|=|(A_1 \cup A_2
\cup \cdots \cup A_{d }) \cap S_1| \leq \frac{d l}{2}.\]
 However, the choice of $A_1, A_2, \ldots, A_d$ ensures that $A_1 \cap A_2
\cdots \cap A_{d } \cap S_2 = \emptyset$, contradicting Condition
(2). This leads to the conclusion that there exists at most one set
$S_m$ such that
\[ | \mathcal{F}\cap [S_{m}]^l| \geq  \frac{d}{2}.\]
Without loss of generality,  let us assume that $m = t$. Thus we
have
\[
|\mathcal{F}\cap [S_i]^l| \leq \frac{d-1}{2},
\]
for $i=1,\ldots ,t-1$. Set \[ \mathcal{H}_i=\{F\in \mathcal{F} :
|F\cap S_i|=l-1\}\]
 and \[ \deg_{\mathcal{H}_i} (B)=|\{F\in
\mathcal{H}_i: B\subset F\}|\]
 for each $1 \leq i \leq t$.

We claim that there exists at least one set $S_i$ ($i\in \{1,\ldots
,t\}$) such that
\[
|\mathcal{H}_i|\leq {k\choose {l-1}} \quad \mbox{and} \quad
|\mathcal{F}\cap [S_i]^l|\leq \frac{d-1}{2}.
\]
Suppose that the above claim is not true. Then
\begin{equation}\label{lem331}
|\mathcal{H}_i| \geq{k\choose {l-1}}+1,
\end{equation}
for $i=1,\cdots ,t-1$. Moreover, if $|\mathcal{F}\cap [S_t]^l|
\leq \frac{d-1}{2}$, then
\[|\mathcal{H}_t| \geq{k\choose {l-1}}+1.\]
By (\ref{lem331}), there exists a $(l-1)$-subset $B$ of $S_1$ such
that
\begin{equation}\label{lem332}
\deg_{\mathcal{H}_{1}}(B)\geq 2.
\end{equation}
Assume that $A_1,A_2\in \mathcal{H}_1$ are chosen subject to the
conditions  $B\subset A_1$ and $B\subset A_2$. Since
\[ |\mathcal{H}_2| \geq {k\choose {l-1}}+1>d-2,\]
 we can choose
$A_3,\ldots A_d$ from $\mathcal{H}_2$. Since $A_1 \cap A_2 = B
\subseteq S_1$,
\[
A_1\cap \cdots \cap A_d\cap S_2=\emptyset
\]
and
\[
|A_1\cup \cdots \cup A_d- S_2|\leq (l+1)+(d-2)=l+d-1 \leq
\frac{dl}{2}
\]
for $d \geq 4$ and $l\geq 3$. So we have reached a contradiction
to Condition (2) when $d \geq 4$.

Consider the case $d=3$. Let $\{x_i\}=A_i-B$ for $i=1,2$. Since
$A_1,A_2\in \mathcal{H}_1$, we have $x_i \not \in S_1$. Let $x_1 \in
S_{i_0}$ for some $i_0\geq 2$. Choose $A_3$ to be either in
$\mathcal{H}_{i_0}$ or $\mathcal{F}\cap [S_{i_0}]^{l}$. Since
$A_1\cap A_2=B\in S_1$ and $S_1 \cap S_2=\emptyset$, we have
\[
A_1\cap A_2 \cap A_3\cap S_{i_0}=\emptyset
\]
and
\[
|A_1\cup A_2 \cup A_3- S_{i_0}|\leq (l-1)+1+1=l+1 \leq
\frac{dl}{2}
\]
for $l\geq 3$ and $d = 3$, contradicting Condition (2) again. Thus
the claim is verified.

Without loss of generality, we assume that
\[
 |\{F\in \mathcal{F} : |F\cap S_1|=l-1\}|=|\mathcal{H}_1| \leq
{k\choose {l-1}} \quad \mbox{and} \quad |F\cap [S_1]^l|\leq
\frac{d-1}{2}.
\]
For any $F \in \mathcal{F}$, we may express $F$ as $F_1 \cup F_2$,
where $F_1 = F \cap S_1$ and $F_2 = F - F_1$. For a fixed $F_1$ of
size $l-r$ ($1 \leq r \leq l$), let $\mathcal{F}_r$ be the family of
all $r$-sets $F_2 \subset S_2 \cup S_3 \cup \cdots \cup S_t$ such
that $F_1 \cup F_2 \in \mathcal{F}$.

We claim that $\mathcal{F}_{r}$ satisfies the conditions of the
lemma. For otherwise, we may assume that there exist $A_1$, $A_2$,
$\dots$, $A_{d }$ $\in \mathcal{F}_{r}$ and $i \in \{2,\cdots
,t\}$ such that $A_1 \cap A_2 \cap \cdots \cap A_{d } \cap S_i =
\emptyset$ and
 \[ |(A_1 \cup A_2 \cup
\cdots \cup A_{d }) - S_i| \leq \frac{d }{2}r.\]
 Now, let
$A_{j}'=A_j \cup F_1 $ for $1\leq j \leq d $. Clearly, $A_1', A_2',
\ldots, A_d' \in \mathcal{F}$ and

\noindent $A_1' \cap A_2' \cap \cdots \cap A_{d }' \cap S_i =
\emptyset$. Recalling that $l\geq r$, we find
\[|(A_1' \cup A_2' \cup \cdots \cup A_{d }') - S_i|=|F_1|+|(A_1 \cup A_2 \cup
\cdots \cup A_{d }) - S_i|\] \[ \leq l-r + \frac{d
r}{2}=l+\frac{d-2}{2}r\leq l + \frac{d-2}{2}l =\frac{d l}{2},\]
contradicting Condition (2). Thus we have shown that
$\mathcal{F}_{r}$ satisfies the conditions of the lemma. For $r
\geq 2$, by the inductive hypothesis, we see that \[
|\mathcal{F}_r|< {{(t-1)k-1} \choose {r-1}}.\] Since $l \geq 3$
and $d  \leq k$, it is easy to check that
\[\sum_{r=2}^{l}{{k} \choose {l-r}} - d \geq 0.\]
Hence $|\mathcal{F}|$ can be bounded as follows,
 \begin{eqnarray*}
|\mathcal{F}| & \leq & \sum_{r=2}^{l} {{k} \choose
{l-r}}|\mathcal{F}_{r}|+|\{F\in \mathcal{F} : |F\cap
S_1|=l-1\}|+|\mathcal{F}\cap[S_1]^{l}|
\\[6pt]
& \leq & \sum_{r=1}^{l} {{k} \choose {l-r}}{{(t-1)k-1} \choose
{r-1}}-\sum_{r=1}^{l}{{k} \choose {l-r}} + {k\choose
{l-1}}+\frac{d-1}{2} \\[6pt]
  & <  & {{tk-1} \choose
{l-1}}-\sum_{r=2}^{l}{{k} \choose {l-r}} + d \leq {{tk-1} \choose
{l-1}}.
\end{eqnarray*}
This completes the proof. \qed

We are now ready to prove Theorem 2.2.

\noindent {\it Proof of Theorem 2.2.} For $d =3$, the result follows
from Theorem 1.8. So we assume $d \geq 4$. Let $S_1$, $S_2$,
$\dots$, $S_t$ be a maximum subfamily of pairwise disjoint
$k$-subsets from $\mathcal{F}$. We proceed by  induction on $t$. If
$t=1$, then $\mathcal{F}$ is intersecting and the result follows
from Theorem 1.1 when $n\geq 2k$. When $\frac{dk}{d-1}\leq n<2k$,
for any $A_1,\ldots ,A_d\in \mathcal{F}$, $|A_1\cup \cdots \cup
A_d|\leq n<2k$, it follows that their intersection is nonempty from
the condition of the theorem. Hence the theorem reduces to Theorem
1.2 in this case. Now we may assume that $t \geq 2$ and the theorem
holds for $t-1$. Note that $t=1$ is the only case when $\mathcal{F}$
can be a complete star. It will be shown that $|\mathcal{F}| <
{{n-1} \choose {k-1}}$.

If $n = tk$,  we set $l = k$. The condition on $\mathcal{F}$ in
Theorem 2.2 implies the conditions on $\mathcal{F}$ in Lemma 2.4
with $d$ replaced by $d-1$. In fact, suppose that there exist $A_1$,
$A_2$, $\dots$, $A_{d-1}$ $\in \mathcal{F}$ for which $A_1 \cap A_2
\cdots \cap A_{d-1} \cap S_i = \emptyset$. Since every $(d,
{d+1\over 2})$-cluster of $\mathcal{F}$ is intersecting,  we see
that
\[ |A_1 \cup A_2 \cup \cdots \cup A_{d-1} \cup
S_i| > \frac{d-1}{2}k,\]
 hence
\[ |A_1 \cup A_2 \cup \cdots \cup A_{d-1}- S_i| > \frac{d+1}{2}k -k
= \frac{d-1}{2}k.\] Hence the theorem follows from Lemma 2.4 in this
case.

We now assume $n > tk$ and let
 \begin{equation}\label{y}
  Y = [n] -
\bigcup_{i=1}^{t}S_i.\end{equation}   Given the choice of $S_1$,
$S_2$, $\dots$, $S_t$, $Y$ does not contain any subset $A \in
\mathcal{F}$. Set \[ \mathcal{F'}=\{F \in \mathcal{F}\colon |F \cap
Y| = k-1\}.\]

We claim that if  $|Y| = n - tk \geq k$, then
\begin{equation}\label{cla}
|\mathcal{F'}|\leq {{n-tk}\choose {k-1}}.
\end{equation}
If the claim is not true given the condition, then we have \[
|\mathcal{F'}|\geq {{n-tk}\choose {k-1}}+1 \geq k+1
> d.\]  Therefore, there exists a $(k-2)$-subset $B \subset Y$ such
that
\begin{equation}\label{deg}
\deg_{\mathcal{F}'}(B)\geq |Y|-k+3 = (n-tk) -k+3.
\end{equation}
Otherwise, we would have
\[
|\mathcal{F'}|\leq \frac{((n-tk)-k+2){{n-tk}\choose
{k-2}}}{k-1}={{n-tk}\choose {k-1}}.\]
 Since the number of
 $(k-1)$-subsets of $Y$ containing $B$ is equal to $|Y|-k+2$, there exists =
$(k-1)$-subset $C$ in $Y$ containing $B$ such that
$\deg_{\mathcal{F}'}(C)\geq 2$. Let $A_1,A_2 \in \mathcal{F}'$ be
such that $A_1 \cap A_2 = C \subset Y$. It is easy to see that
\[
A_1\cap A_{2}\cap S_i=\emptyset
\]
for each $1 \leq i \leq t$. Let $A_3,A_4, \dots, A_{d-1}$ be
additional subsets in $ \mathcal{F}'$ such that $B \subseteq A_i$
for each $i$ if $|Y|-k+3\geq d-1$. We deduce that
\[
A_1\cap \cdots \cap A_{d-1}\cap S_i=\emptyset
\]
for each $1 \leq i \leq t$. Moreover,
 \[  |A_1\cup \cdots \cup
A_{d-1}|\leq k-2+2(d-2)+1=k+2d-5,\;\; \mbox{if } |Y|-k+3\geq d-1
\] and \[ |A_1\cup \cdots \cup A_{d-1}|\leq |Y|+d-1 \leq k+2d-6,
\;\; \mbox{if } |Y|-k+3 < d-1.
\]
Let $S_h$ be such that $S_h \cap A_1 \neq \emptyset$. Since $k\geq d
\geq 4$, we see that
\[
|(A_1\cup \cdots \cup A_{d-1}) \cup S_h|  \leq k+2d-5 +
(k-1)=2k+2d-6 \leq \frac{d+1}{2}k,
\]
contradicting the assumption of the theorem. So  the claim is
justified.

Note that for any member $F$  in $\mathcal{F}$, we can  write it as
$F=F_1 \cup F_2$, where $F_1 = F \cap Y$ and $F_2 = F - F_1$. We now
consider all possible ways to construct $F$ in the above form. Let
$F_1$ be a given subset of $Y$ size $k-l$ ($1 \leq l \leq k$). By
the definition of $Y$ in (\ref{y}), $F_2$ is a subset
$\cup_{i=1}^{t}S_i$. Let $\mathcal{F}_l$ be the family of all
$l$-sets $F_2 \subset \cup_{i=1}^{t}S_i$ such that $F_1 \cup F_2 \in
\mathcal{F}$. It remains to prove that $\mathcal{F}_l$ satisfies the
conditions in Lemma 2.4 with $d$ replaced by $d-1$. For $l=k$, the
assumption of the theorem implies that for every $A_1, A_2, \dots,
A_{d-1} \in \mathcal{F}_{k}$, if $A_1 \cap A_2 \cap \cdots \cap
A_{d-1} \cap S_i = \emptyset$, then
\[ |A_1 \cup A_2 \cup \cdots \cup A_{d-1} \cup S_i|
> \frac{d+1}{2}k\]
 which yields that
 \[ |A_1 \cup A_2 \cup \cdots \cup
A_{d-1} - S_i| > \frac{d-1}{2}k.\] Therefore, the assertion holds
when $l=k$. For $l < k$, if the assertion is not valid, then there
exist $A_1, A_2, \dots, A_{d-1} \in \mathcal{F}_{l}$ such that $A_1
\cap A_2 \cdots \cap A_{d-1} \cap S_i = \emptyset$ and
\[ |A_1 \cup A_2 \cup \cdots \cup A_{d-1} - S_i| \leq \frac{d-1}{2}l.\]
Setting $A_{i}'=A_i \cup F_1 $ for $i \leq d-1$, we deduce that
$A_i' \in \mathcal{F}$, $A_1' \cap A_2' \cdots \cap A_{d-1}' \cap
S_i = \emptyset$, and
\[|(A_1' \cup A_2' \cup \cdots \cup A_{d-1}') \cup S_i|=|F_1|+|(A_1 \cup =_2 \cup
\cdots \cup A_{d-1}) - S_i| + |S_i|\]
\[ \leq k-l +
\frac{d-1}{2}l+k = 2k + \frac{d-3}{2}l\leq 2k + \frac{d-3}{2}k
=\frac{d+1}{2}k,\] contradicting the assumption of the theorem. Up
to now, we have shown that $\mathcal{F}_{l}$ satisfies the
conditions in Lemma 2.4. For $l \geq 2$, by Lemma 2.4 we find that
\[ |\mathcal{F}_{l}|< {{tk-1} \choose {l-1}}.\]
Evidently, for $|Y|=n-tk \leq k-2$, we have \[ |\{F \in
\mathcal{F}\colon |F \cap Y| = k-1\}|=0.\]
 For the case $|Y| = k-1$, we have
 \[ |\{F \in \mathcal{F}\colon |F \cap Y| =
k-1\}|< d-1 \leq k-1.\] Otherwise we can choose $d-1$ sets
$A_1,\ldots ,A_{d-1} \in \mathcal{F}$ together with $S_1$ in
violation of the assumption of theorem. When $|Y| \geq k$, It
follows from (\ref{cla}) that
\[|\{F \in \mathcal{F}\colon |F \cap Y| = k-1\}| \leq {{n-tk}\choose
{k-1}},\] which implies
\[|\{F\in \mathcal{F}\colon |F \cap Y| = k-1\}| < \sum_{l=1}^{k} {{n-tk} \choose
{k-l}}.\] Finally,
\begin{eqnarray*}
 |\mathcal{F}| & \leq  & \sum_{l=2}^{k} {{|Y|} \choose
{k-l}}|\mathcal{F}_{l}|+|\{F \in \mathcal{F}: |F \cap Y| = k-1\}|\\[6pt]
&  \leq  & \sum_{l=2}^{k} {{|Y|} \choose {k-l}}\left[{{tk-1}
\choose {l-1}}-1\right] + |\{F \in \mathcal{F}: |F \cap Y| =
k-1\}|\\[6pt]
& =  & \sum_{l=1}^{k} {{|Y|} \choose {k-l}}\left[{{tk-1}
\choose {l-1}}-1\right] + |\{F \in \mathcal{F}: |F \cap Y| =
k-1\}|\\[6pt]
& = & \sum_{l=1}^{k} {{n-tk} \choose {k-l}}{{tk-1} \choose {l-1}}
- \sum_{l=1}^{k} {{n-tk} \choose {k-l}}+|\{F \in \mathcal{F}: |F
\cap
Y| = k-1\}| \\[6pt]
&<& {{n-1} \choose {k-1}},
\end{eqnarray*}  as required. This
completes the proof.\qed

\vspace{2mm}

\noindent {\bf Acknowledgments.} The authors  wish to thank the
referees for their helpful suggestions. This work was supported by
the 973 Project, the PCSIRT Project of the Ministry of Education,
the Ministry of Science and Technology, and the National Science
Foundation of China.

\end{document}